\documentclass[12pt]{amsart}
\usepackage{epsfig,psfrag,verbatim,amssymb,amsfonts}
\newtheorem{theorem}{Theorem}
\newtheorem{lemma}[theorem]{Lemma}
\newtheorem{prop}[theorem]{Proposition}
\newtheorem{cor}[theorem]{Corollary}

\newtheorem{ex}{Example}
\newtheorem{rem}{Remark}
\newtheorem*{problem}{Open Problem}

\newcommand{\field}[1]{\mathbb{#1}}
\newcommand{\R}{\field{R}}

\newcommand{\Z}{\field{Z}}
\newcommand{\N}{\field{N}}
\newcommand{\PP}{\field{P}}

\newcommand{\EE}{\field{E}}

\newcommand{\D}{D}

\newcommand{\F}{\mathcal{F}}

\newcommand{\si}{\sigma}
\newcommand{\al}{\alpha}
\newcommand{\ga}{\gamma}
\newcommand{\om}{\omega}
\newcommand{\eps}{\varepsilon}
 \newcommand{\won}{{\boldsymbol 1}}
\newcommand{\foot}{}

\newcounter{constante}
\newcommand{\con}[1]{
                    \immediate\write 1{\noexpand\newlabel{#1}{{\theconstante}{\theconstante}}}
                    c_{\theconstante}
                    \stepcounter{constante}
                   }

\newcommand{\abel}[1]{}

\begin{document}

\setcounter{page}{1}

\title[Multi-excited random walks]
{
Multi-excited random walks on integers}

\maketitle

\begin{center}
{\sc By Martin P.W.\ Zerner
}

\end{center}\vspace*{5mm}

\begin{center}\begin{quote}
{\footnotesize {\sc Abstract}. We introduce a class of nearest-neighbor
integer random walks in random and non-random
media, which 
includes excited random walks considered in the literature.
At each site the random walker has a drift to the right,
the strength of which depends on the environment at that site and on
how often the walker has visited that site before.
We give exact criteria for
recurrence and transience and consider the speed of the walk. 
}\vspace*{5mm}
\end{quote}
\end{center}

\section{Introduction}
\begin{figure}[t]\label{robin}
\epsfig{file=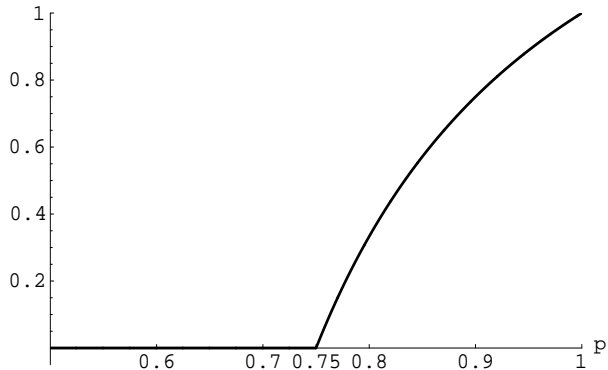,height=5cm, width=8cm, angle=0}
 \caption{\foot The probability that a random walker will never return to its starting point if there are two cookies
with parameter $p$ each at each integer.}
\end{figure}

The results of the present paper are best illustrated by the following example.
\begin{ex}\label{zwei}
{\rm We put two cookies on each integer and launch a nearest neighbor random walker at the origin.
Whenever there is at least one cookie at the random walker's current position, the walker eats exactly one of these 
cookies, thus removing it from this site, and then jumps independently of its past to
the right with probability $p$ and to the left with probability $1-p$, where $p\in[1/2,1]$ is a fixed parameter.
Whenever there is no cookie left at the random walker's current position, the walker jumps 
 independently of its past to the
left or right with equal probability 1/2. 

We shall show a phase transition in the recurrence and transience behavior of the walk, see Theorem \ref{main}: 
 If $1/2\leq p\leq 3/4$ then
the walker will visit its starting point 0
almost surely infinitely often. However, if $p>3/4$ then the walker will 
visit 0 almost surely only finitely many times. Moreover, the probability that the walker will never return 
to 0 is $(1-2/(2p-1))_+$, see Figure \ref{robin} and Theorem \ref{hp}. 
Finally,
for all $p<1$ the walk has zero speed, even if it is transient,
see Theorem \ref{v0}. 
}\hfill $\Box$
\end{ex}\abel{zwei}
\abel{robin}
This example can be formalized and generalized as follows.
A cookie environment is an element 
\[\om=(\om(z))_{z\in\Z}=\left((\om(z,i))_{i\geq 1}\right)_{z\in\Z}\quad\in\quad
\Omega_+:=\left([1/2,1]^{\N}\right)^{\Z}.\]
We will refer to $\om(z,i)$ as to the strength of the $i$-th cookie at $z$. \
This is the probability for the
random walker to jump from $z$ to $z+1$ if it is currently visiting $z$ for the $i$-th time.
T.\ Komorowski suggested to consider the  cookies as \textit{bribes} which push an otherwise
unbiased walker to the right. 

More formally, given a starting point $x\in\Z$ and a cookie environment $\om\in\Omega_+$, we consider
an integer valued 
 process $(X_n)_{n\geq 0}$ on some suitable probability space $\left(\Omega, \F, P_{x,\om}\right)$ 
for which  the process of its history $(H_n)_{n\geq 0}$
defined by
$H_n:=(X_m)_{0\leq m\leq n}$ is a Markov chain, which satisfies
$P_{x,\om}$-a.s.\ 
\begin{eqnarray*}
P_{x,\om}[X_0=x]&=&1,\\
P_{x,\om}[X_{n+1}=X_n+1\mid H_n]&=&\om\left(X_n,\#\{m\leq n\mid X_m=X_n\}\right),\\
P_{x,\om}[X_{n+1}=X_n-1\mid H_n]&=&1-\om\left(X_n,\#\{m\leq n\mid X_m=X_n\}\right).
\end{eqnarray*}
Note that $(X_n)_n$ itself is in general not a Markov chain since its transition probabilities depend on
the history of the process. 
In Example \ref{zwei} we have chosen $x=0$ as starting point and the cookie environment 
$\om\in\Omega_+$ with
$\om(z)=(p,p,1/2,1/2,1/2,\ldots)$ for all $z\in\Z$.

This model generalizes in part one-dimensional excited random walks (ERW) and
random walk perturbed at its extrema, also called $pq$ walks,  see Benjamini-Wilson \cite{bewi} and 
Davis \cite{dav1}
for results and references, regarding also continuous space and time analogues.
For higher dimensional ERWs see \cite{bewi},
  Kozma \cite{kozma}, and Volkov \cite{vol}.

 The intersection of our model, which we call multi-ERW, with one-\-dimen\-sional 
ERW as defined in \cite {bewi} and \cite{dav1} 
deals, in our language,
with cookie environments of the form $\om(z)=(p,1/2,1/2,1/2,\ldots)$ for all $z\in\Z$,
where $p\in[1/2,1]$ is fixed. In such an environment the walker is excited, i.e. biased
to the right,  only on the first
visit to a site.
We call such random walks once-excited.

The novelty of our model of multi-ERW  is that it permits  
different levels of excitement for different visits 
to a site. Moreover, the excitement levels may vary randomly from site to site, see Section \ref{ratire}
for details. While once-ERW is recurrent for all $p<1$  (see \cite[p.\ 86]{bewi}),
multi-ERW exhibits a more interesting recurrence and transience behavior,
as highlighted in Example \ref{zwei}.

Our motivation for the study of multi-ERW on integers came from 
the problem posed at the end of \cite{bewi}, as to whether once-ERW on
$\Z^2$ has positive speed. Although we do not see how to prove this,
morally, this should be true if once-ERW
on a strip $\Z\times\{0,\ldots,k\}$ has finite speed for $k$ large enough. Moreover, once-ERW
 on a strip of finite width should roughly behave like 
multi-ERW on $\Z$ with a finite number of cookies per site. 

A second source of motivation was to find a unifying model which includes both
once-ERW and random walks in random environments (RWRE, see e.g.\ \cite{so}, \cite{sz}, \cite{zei}) as special cases,
see Remark \ref{rem2} in Section \ref{ratire} for details.

Let us now describe how the remainder of the present paper is organized.
Section \ref{notation} provides basic lemmas which will be used throughout the paper. 
After some preparation we will also describe in Remark \ref{rem1} the main idea behind the proof of the phase transition
described in Example \ref{zwei}. 
In Section  \ref{ratide} we introduce the notion of recurrence and transience of
states in fixed environments $\om$. This will be used in Section \ref{ratire},
which contains our main result Theorem \ref{main}, a sufficient and necessary criterion 
for recurrence in stationary and ergodic environments. In Section \ref{elo} we investigate
 random walks which one after the other live on the environment
left over by the previous random walk. Section \ref{sloln} is devoted to a strong
law of large numbers 
for the walk in a stationary and ergodic environment.  
 Section \ref{m} deals
with the monotonicity of the return probability 
and the speed,
two quantities, which are explicitly computed in the last
section for the case in which the excitement is gone after the second visit.

\section{Notation and Preliminaries}\label{notation}\abel{notation}
Let
\[T_k:=\inf\{n\geq 0\mid X_n=k\}\]
be the first passage time of $k$. 
The following lemma will be generalized in Lemma \ref{monoenv} by a different technique.
\begin{lemma}\label{simply}
For all $x<y<z$ and all $\om\in\Omega_+$,
\[
P_{y,\om}[T_x<T_z]\leq \frac{z-y}{z-x}.
\]
\end{lemma}\abel{simply}
In particular, by letting $x\to-\infty$ we see that $T_z$ is $P_{y,\om}$-a.s.\ finite.
\begin{proof}
We couple $(X_n)_{n\geq 0}$ to a simple symmetric random walk $(Y_n)_{n\geq 0}$ starting at $y$ such that
almost surely $Y_n\leq X_n$ for all $n\geq 0$. To this end, we may assume 
that there is a sequence  $(U_n)_{n\geq 0}$ of independent random variables on $\Omega$ which are
uniformly distributed on $[0,1]$.
If  the walk $(X_n)_n$ visits at time $n$ a site $x$ for the $j$-time $(j\geq 1)$ then it moves to the
right in the next step iff $U_n<\om(x,j)$, whereas the walk $(Y_n)_n$ jumps to the right iff $U_n<1/2$.
Then $(X_n)_n$ is an ERW in the environment $\om$ whereas $(Y_n)_n$ is a simple symmetric
random walk.
Since $\om(x,j)\geq 1/2$ we get $Y_n\leq X_n$ almost surely by induction over $n$.
Therefore, if $(X_n)_n$ exits the interval $]x,z[$ in $x$ then so does $(Y_n)_n$, which has probability $(z-y)/(z-x)$.
\end{proof}
The average displacement of the walk after having eaten a cookie of strength $p$ is $2p-1$. Therefore, 
\[\delta^{x}(\om):=\sum_{i\geq 1}(2\om(x,i)-1)\]
is the total drift stored in the cookies at site $x$ in the environment $\om$.
The drift contained in the cookies at site $x$ which have been eaten before time $n$ will
be called
\[\D_{n}^x:=\sum_{i=1}^{\#\{m<n\mid X_m=x\}}(2\om(x,i)-1).\]
To distinguish between recurrence and transience we will distinguish between cookies on
nonnegative and negative integers. Therefore, we introduce
\[
\D_n^+:=\sum_{x\geq 0}\D_{n}^x,\qquad
\D_n^-:=\sum_{x< 0}\D_{n}^x,\quad\mbox{and}\quad \D_n :=\D_n^++\D_n^-.
\]
\begin{lemma}\label{expect}Let  $\om\in\Omega_+$ such that
\begin{equation}\label{pine}
\liminf_{i\to\infty}\frac{1}{i}\sum_{y=-i}^{0}(2\om(y,1)-1)>0.
\end{equation}\abel{pine}
Then for all $x,k\in\Z$ with $k\geq x$,
\begin{equation}\label{hallo}
E_{x,\om}[\D_{T_k}]=k-x.
\end{equation}\abel{hallo}
\end{lemma}\abel{expect}
Note that for simple symmetric random walk, i.e.\ for $\om\equiv 1/2$, $\D_{T_k}=0$\ $P_{x,\om}$-a.s..
Hence assumption (\ref{pine}) is essential.
\begin{proof}
By shifting $\om$ by $x$ to the left we may assume 
 without loss of generality $k\geq x=0$.
Consider  the process 
$M_n:=X_n-\D_n$\ $(n\geq 0)$.
It is standard to check that $(M_n)_{n\geq 0}$ is 
a $P_{0,\om}$-martingale  with respect to the filtration  $(\F_{n})_{n\geq 0}$
generated by $(X_n)_{n\geq 0}$. 
Therefore,
 by  the Optional Stopping Theorem for all $n\geq 0$,
\[0=E_{0,\om}[M_{T_k\wedge n}]=E_{0,\om}[X_{T_k\wedge n}, T_k\leq n]
+E_{0,\om}[X_{T_k\wedge n}, n<T_k]-E_{0,\om}[\D_{T_k\wedge n}]
\]
 and consequently,
\begin{equation}\label{hai}
E_{0,\om}[\D_{T_k\wedge n}]=
kP_{0,\om}[T_k\leq n]+E_{0,\om}[X_{n}, n<T_k].
\end{equation}\abel{hai}
Now consider (\ref{hai}) as $n\to\infty$.
Since $\om\in\Omega_+$, 
the left hand side of (\ref{hai}) tends  by monotone convergence to $E_{0,\om}[\D_{T_k}]$.
Moreover, the first term on the right-hand side goes to $k$. 
Consequently,
\begin{equation}\label{hawaii}
E_{0,\om}[\D_{T_k}]=
k+\lim_{n\to\infty}E_{0,\om}[X_{n}, n<T_k].
\end{equation}\abel{hawaii}
Hence all that remains to be shown is 
\begin{equation}\label{all}
\lim_{n\to\infty}E_{0,\om}[X_{n}, n<T_k]=0.
\end{equation}\abel{all}
Since for all $n\geq 0$,
\begin{equation}\label{biele}
\min_{m<T_k}X_m\leq X_{n}\won\{n<T_k\}\leq k\qquad\mbox{ $P_{0,\om}$-a.s.,}
\end{equation}\abel{biele}
(\ref{all}) will follow by  dominated convergence once we have shown that the non-negative
random variable $-\min_{m<T_k}X_m$ has finite $E_{0,\om}$-expectation.
Denote by $\ga$ the left hand side of (\ref{pine}). Then
\begin{eqnarray}\label{irrigation}
E_{0,\om}\Big[-\min_{m<T_k}X_m\Big]&\leq& E_{0,\om}[2\D_{T_k}/\ga]\\
&&+\ E_{0,\om}\Big[-\min_{m<T_k}X_m,\ 
-\min_{m<T_k}X_m>2\D_{T_k}/\ga\Big].\label{ditch}
\end{eqnarray}\abel{irrigation}\abel{ditch}
The term on the right hand side of (\ref{irrigation}) is finite since 
  (\ref{hawaii}) and (\ref{biele}) imply 
$E_{0,\om}[\D_{T_k}]\leq k.$
The term in (\ref{ditch}) equals
\begin{equation}\label{only}
\sum_{i\geq 1}i P_{0,\om}\Big[-\min_{m<T_k}X_m=i, \D_{T_k}<\ga i/2\Big].
\end{equation}\abel{only}
Observe that on the event 
$\{T_{-i}<T_k\},$
\[\D_{T_k}\geq \sum_{y=-i}^{0}(2\om(y,1)-1).\]
Therefore, (\ref{only}) is less than or equal to 
\begin{equation}\label{gewitter}
\sum_{i\geq 1}i \won\left\{\frac{1}{i}\sum_{y=-i}^{0}(2\om(y,1)-1)< \frac{\ga}{2} \right\},
\end{equation}\abel{gewitter}
which is finite since due to the choice of $\ga$  
only finitely many indicator functions in (\ref{gewitter}) do not vanish.
\end{proof}

\begin{rem}\label{rem1}
{\rm We are now ready to present the idea of the proof of the recurrence and transience behavior
in the two-cookie case described in Example \ref{zwei}. This will be made rigorous and more
general in Theorem \ref{main}. Roughly speaking, (\ref{hallo})  states that an ERW
starting at 0 needs to eat $k/(2p-1)$ cookies in order to reach $k$. Compare this number to the total number $2k$
of cookies available between 0 and $k-1$.
If $2k<k/(2p-1)$ then 
the walker needs to visit once in a while negative integers in order to meet its cookie needs because
there are not enough cookies available on the positive integers. This makes the
walker recurrent.

On the other hand, if $2k>k/(2p-1)$ then the walker cannot afford to return to 0 infinitely often
because on its way back from its up-to-date maximum value, say $k-1$, to 0 the ERW will eat all the remaining cookies between
0 and $k-1$, thus consuming at least $2k$ cookies before it reaches $k$. This would be more than the $k/(2p-1)$ 
cookies the ERW
should eat. Therefore, the walker has to
be transient.}\hfill $\Box$
\end{rem}

In the following we are concerned with the probabilities of the events
\[R_k:=\{X_n=k\ i.o.\}=
\limsup_{n\to\infty}\{X_n=k\}\quad(k\in\Z)\]
that any given site $k$ is visited infinitely often.
The following lemma 
states  that the behavior of the walk to the right of $k$ does not depend on the environment
to the left of $k$ nor on where to the left of $k$ the walk started. 
 A related result for Brownian motion
perturbed at its extrema
is \cite[Proposition 1]{perwer}. 
Consider the sequences $(\tau_{k,m})_{m\geq 1}$\ $(k\in\Z)$
defined by 
\[\tau_{k,0}:=-1\quad\mbox{ and}\quad 
\tau_{k,m+1}:=\inf\{n>\tau_{k,m}\mid X_n\geq k\}.
\] 
They enumerate the times $n$ at which $X_n\geq k$. Note
that these times are  stopping times with respect to $(\F_n)_{n\geq 0}$.
Moreover, they are $P_{x,\om}$-a.s.\ finite $(x\in\Z)$ since $T_r<\infty$ for all $r\geq 0$.
\begin{lemma}\label{same}
Let $x_1,x_2\leq k$ and  $\om_1,\om_2\in\Omega_+$ such that $\om_1(x)=\om_2(x)$ 
for all $x\geq k$.
Then
$(X_{\tau_{k,m}})_{m\geq 0}$ has the same distribution under $P_{x_1,\om_1}$ as under 
 $P_{x_2,\om_2}$. 
In particular, 
\begin{equation}\label{wolfgang}
P_{x_1,\om_1}[R_k]=P_{x_2,\om_2}[R_k].
\end{equation}
\end{lemma}\abel{same}\abel{wolfgang}
In the proof of Lemma \ref{same}  and throughout the paper we will use the (strong) Markov property
for the Markov chain $(H_n)_n$. To this end we need to introduce notation for the
cookie environment left behind by a cookie eating random walker.
For any $\om\in\Omega_+$ and any finite  sequence $(x_n)_{n\leq m}$ of integers we define 
$
\psi\left(\om, (x_n)_{n\leq m}\right)\in\Omega_+
$
by
\begin{equation}\label{psi}
\psi(\om,(x_n)_{n\leq m})(x,i):=\om\left(x,i+\#\{n<m\mid x_n=x\}\right).
\end{equation}\abel{psi}
This is the environment we obtain form $\om$ by following the path $(x_n)_{n<m}$ and removing
the bottom cookie in each site visited. Note that in the definition
of $\psi$ we do not remove the cookie from the final site $x_m$.

\begin{proof}[Proof of Lemma \ref{same}]
It suffices to show that for all sequences $(y_m)_{m\geq 1}$ with $y_m\geq k\ (m\geq 1)$ and for all $M\geq 1$,
\begin{equation}\label{hawk}
P_{x_1,\om_1}[A_M]=P_{x_2,\om_2}[A_M],
\end{equation}\abel{hawk}
 where 
\[A_M:=A_M\left((y_m)_{m\geq 1}\right):= \left\{\left(X_{\tau_{k,m}}\right)_{m=1}^M=(y_m)_{m=1}^M\right\}
\quad(M\geq 1).\]
So fix such a sequence $(y_m)_{m\geq 1}$.
For $M=1$, (\ref{hawk}) is trivial since $X_{\tau_{k,1}}=k$ $P_{x_i,\om_i}$-a.s.\ $(i=1,2)$. 
Now assume that (\ref{hawk}) has been proven for $M$. Then
\begin{eqnarray}
P_{x_i,\om_i}[A_{M+1}]&=&
E_{x_i,\om_i}\left[P_{x_i,\om_i}\left[X_{\tau_{k,M+1}}=y_{M+1}\mid \F_{\tau_{k,M}}\right],
A_M\right]\nonumber \\ 
&=&\label{va}
E_{x_i,\om_i}\left[C, A_M\right],
\end{eqnarray}\abel{va}
where  by the strong Markov property 
\[C:=P_{y_M,\psi\left(\om_i,H_{\tau_{k,M}}\right)}\left[X_{\tau_{k,2}}=y_{M+1}\right].\]
It suffices to show  that on  $A_M$, $C$ is equal to a deterministic constant $c$ which may depend only on 
$(y_m)_{m\geq 0}$ and on $\om_1$ and $\om_2$ where they coincide, thus being  independent of $i$. 
Indeed, then the right-hand side of (\ref{va})
is equal to $cP_{x_i,\om_i}[A_M]$, which is independent of $i$ by induction hypothesis.
Since $(X_n)_n$ is a nearest neighbor walk, $C=0$ unless $y_{M+1}\in\{y_M+1, (y_M-1)\vee k\}$.
If  $y_{M+1}= (y_M-1)\vee k$ then on $A_M$,
\begin{eqnarray*}
C&=& 1-\psi\left(\om_i,H_{\tau_{k,M}}\right)(y_M,1)\\
&=&1-\om_i\left(y_M,1+\#\left
\{n<\tau_{k,M}\mid X_n=y_M\right\}\right)\\
&=&1-\om_i\left(y_M,1+\#\left
\{m<M\mid y_m=y_M\right\}\right).
\end{eqnarray*}
 Since by assumption $\om_1(y_M)=\om_2(y_M)$, $C$  is independent
of $i$ indeed. A similar argument settles the case $y_{M+1}= y_M+1$.
\end{proof}

\section{Recurrence and Transience in Deterministic Environments} \label{ratide}\abel{ratide}
In this section we establish recurrence and transience criteria for fixed environments.
\begin{lemma}\label{fliege}
Let $x,y,z\in\Z$ with $y< z$ and $\om\in\Omega_+$. Then $P_{x,\om}$-a.s. $R_y\subseteq R_z.$
\end{lemma}\abel{fliege}
\begin{proof} 
Denote by  $T_{y,r}\ (r\geq 1)$ the time of the $r$-th visit to $y$. On the event $R_y$ all times $T_{y,r}\ (r\geq 1)$ 
are finite. Consider the events 
\[ B_{r}:=\{T_z\circ\theta_{T_{y,r-1}}+T_{y,r-1}<T_{y,r}\}\quad(r\geq 0)\]
that the walk visits $z$ between the $(r-1)$-th and the $r$-th visit to $y$. 
Here $\theta_n$ denotes the canonical shift by $n$ steps on the path space.
Since $\{B_r\ \mbox{i.o.}\}\subseteq R_z$ and $B_r\in \F_{T_{y,r}}$ 
it suffices to show by the second Borel Cantelli lemma  (e.g.\ \cite[Ch.4 (3.2)]{durr}) that
$P_{x,\om}$-a.s.\
$\sum_r P_{x,\om}[B_{r+1}\mid \F_{T_{y,r}}]=\infty$ on $R_y$.
To this end, we use the strong Markov property which implies that 
on $R_y$, 
\begin{eqnarray*}
P_{x,\om}[B_{r+1}\mid \F_{T_{y,r}}]&=&P_{y,\psi\left(\om,H_{T_{y,r}}\right)}[T_z<T_{y,2}]\\
&= &\om(y,r)
P_{y+1,\psi\left(\om,H_{T_{y,r}+1}\right)}[T_z<T_{y}]\ \geq\ \frac{1} {2(z-y)}
\end{eqnarray*}
due to Lemma \ref{simply}. 
This is independent of $r$ and hence
not summable in $r$.
\end{proof}
\begin{prop}\label{rectrans}
Let $\om\in\Omega_+$ and $y\in\Z$. Then either for all $x\in\Z$,
$P_{x,\om}[R_y]=0$ or for all $x\in\Z$,
$P_{x,\om}[R_y]=1$.
\end{prop}\abel{rectrans}
If $P_{x,\om}[R_y]=1$  for all $x\in\Z$ then we shall call 
$y$ \textit{$\om$-recurrent}.
Otherwise, i.e.\ if 
$P_{x,\om}[R_y]=0$ for all $x\in\Z$,  $y$ is called
 \textit{$\om$-transient}.
\begin{proof}[Proof of Proposition \ref{rectrans}] 
Let $x\in\Z$ with $P_{x,\om}[R_y]>0$.  All we  have to show is that
\begin{equation}\label{baeh}
\forall k\in\Z\quad P_{k,\om}[R_y]=1.
\end{equation}\abel{baeh} 
By Lemma \ref{fliege}, for all $z> y$,
\[0<P_{x,\om}[R_y]=P_{x,\om}[R_y\cap R_z]\leq P_{x,\om}[(X_n,X_{n+1})=(z,z-1)\ \mbox{i.o.}].\]
Therefore, by the convergence part of the Borel Cantelli lemma,
$\sum_i(1-\om(z,i))=\infty$ for all $z>y$.  However, since the decisions to jump from $z$ to
$z-1$ are made independently of each other under $P_{k,\om}\ (k\in\Z)$, the divergence
part of the Borel Cantelli lemma then implies that for all $k\in\Z$ and all $z> y$
we have $P_{k,\om}$-a.s.\  $R_{z}\subseteq R_{z-1}$.
Since the opposite inclusion holds anyway due to Lemma \ref{fliege}, we have
\begin{equation}\label{lara}
\forall k\in\Z\quad \forall z\geq y\quad R_z\stackrel{P_{k,\om}}{=}R_y.
\end{equation}\abel{lara}
Since $R_y\in\sigma\left(\bigcup_{z\geq 0}\F_{T_{z}}\right)$ the martingale convergence theorem
yields that $P_{x,\om}$-a.s.,
\[
\won_{R_y}=\lim_{z\to\infty}P_{x,\om}\left[R_y\mid\F_{T_z}\right]\ \stackrel{(\ref{lara})}{=}
\lim_{z\to\infty}P_{x,\om}\left[R_z\mid\F_{T_z}\right].
\]
By the strong Markov property this is for all $k\in\Z$ equal to
\[\lim_{z\to\infty}P_{z,\psi\left(\om,H_{T_z}\right)}\left[R_z\right]\  
\stackrel{(\ref{wolfgang})}{=}\ \lim_{z\to\infty}P_{k,\om}\left[R_z\right]
\ \stackrel{(\ref{lara})}{=}\ \lim_{z\to\infty}P_{k,\om}\left[R_y\right]\ =  P_{k,\om}\left[R_y\right],
\]
which implies (\ref{baeh}) because $P_{x,\om}[R_y]>0$ by assumption.
\end{proof}

\begin{ex}{\rm Let $x\in\Z$ and define $\om\in\Omega_+$ by 
 $\om(y,i)=1/2$ if $y\ne 0$ and $\om(0,i)=1-(i+1)^{-2}$ for all $i\geq 1$. Since $\sum_i (i+1)^{-2}$
converges, the Borel Cantelli lemma implies that 
 negative integers are $\om$-transient. On the other hand, nonnegative integers are 
$\om$-recurrent because  simple symmetric random walk is recurrent.
}\hfill $\Box$
\end{ex}

\begin{lemma}\label{eins}
Let $\om\in\Omega_+$ such that 0 is  $\om$-transient.
Then
\[\lim_{K\to\infty}\frac{E_{0,\om}\left[\D_{T_K}^+\right]}{K}=1.
\]
\end{lemma}\abel{eins}
\begin{proof} Since $\won_{R_0}$ and $\D_{T_K}^+$ 
are functions of $\left(X_{\tau_{0,m}}\right)_{m\geq 0}$ and $(\om(x))_{x\geq 0}$
we may change $\om$ due to Lemma \ref{same}
at negative sites without changing  $\om$-transience of $0$ and 
$E_{0,\om}\left[\D_{T_K}^+\right]$. Hence we may 
assume without loss of generality 
that $\om$ satisfies (\ref{pine}).
For $k\geq 1$ consider the possibly infinite stopping time 
\[\si_k:=T_0\circ\theta_{T_{k-1}}+T_{k-1}
\quad\mbox{and the event}\quad A_k:=\{\si_k<T_k\}\]
that the walk  after hitting $k-1$  for the first time, returns to 0 before it reaches $k$.
Note that $A_k\in\F_{T_k}$.
Since 0 is  $P_{0,\om}$-a.s.\ transient, $A_k$ occurs  $P_{0,\om}$-a.s.\ only for finitely
many $k$'s. Hence by the second Borel Cantelli lemma,
\begin{equation}\label{eps}
\sum_{k\geq 1}P_{0,\om}\left[A_k\mid\F_{T_{k-1}}\right]<\infty\quad  \mbox{$P_{0,\om}$-a.s..}
\end{equation}\abel{eps}
Now let $\eps>0$.
Omitting in (\ref{eps}) those $k$'s which are not elements of the set
$S_\eps:=\{k\geq 1\mid  
Y_k > \eps/k\}$, 
where $Y_k:=P_{0,\om}\left[A_k\mid\F_{T_{k-1}}\right],$
we obtain
\[\sum_{k\in S_\eps}\frac{1}{k}<\infty\quad  \mbox{$P_{0,\om}$-a.s..}
\]
Consequently,  $S_\eps$ has $P_{0,\om}$-a.s.\ upper density 0, i.e.\ $\#(\{1,\ldots,K\}\cap S_\eps)/K\to 0$
as $K\to\infty$. Therefore, by dominated convergence,
\begin{equation}
0=\lim_{K\to\infty}E_{0,\om}\left[\frac{\#(\{1,\ldots,K\}\cap S_\eps)}{K}\right]\
=\ \lim_{K\to\infty}\frac{1}{K}\sum_{k=1}^K P_{0,\om}\left[Y_k>\eps/k\right].
\label{brace}
\end{equation}\abel{brace}
 By the strong Markov property, $P_{0,\om}$-a.s.,
\[
Y_k=P_{k-1,\psi\left(\om,H_{T_{k-1}}\right)}[T_0<T_k]\ \leq\ 
\frac{1}{k}
\]
due to Lemma \ref{simply}.
Therefore,
\[E_{0,\om}[Y_k]=E_{0,\om}[Y_k,\ Y_k>\eps/k]+E_{0,\om}[Y_k,\ Y_k\leq\eps/k]\leq \frac{1}{k}P_{0,\om}[ Y_k>\eps/k]+
\frac{\eps}{k}\]
and hence
\[P_{0,\om}[ Y_k>\eps/k]\geq kE_{0,\om}[Y_k]-\eps.\]
Substituting this into (\ref{brace}) yields
\[0\geq\limsup_{K\to\infty}\frac{1}{K}\sum_{k=1}^K(kE_{0,\om}[Y_k]-\eps)\ =\ 
-\eps+\limsup_{K\to\infty}\frac{1}{K}\sum_{k=1}^K kP_{0,\om}[A_k].
\]
Letting $\eps\searrow 0$ gives
\begin{equation}\label{tuch}
0=\lim_{K\to\infty}\frac{1}{K}\sum_{k=1}^K kP_{0,\om}[A_k].
\end{equation}\abel{tuch}
For abbreviation set $\Delta_k^-:=\D_{T_k}^--\D_{T_{k-1}}^-$ for $k\geq 1$. 
This is the total drift of the cookies 
on negative sites which have been eaten between $T_{k-1}$ and $T_k$.
 Since $\Delta_k^-=0$ on $A_k^c$ and  since $A_k\in\F_{\si_k}$ we have
\begin{equation}
\label{ber}
E_{0,\om}[\Delta_k^-]=E_{0,\om}[\Delta_k^-,A_k]=E_{0,\om}\left[E_{0,\om}[\Delta_k^-\mid \F_{\si_k}],A_k\right].
\end{equation}
\abel{ber}
By the strong Markov property this is equal to
\begin{equation}
\label{berg}
E_{0,\om}\left[E_{0,\psi\left(\om,H_{\sigma_k}\right)}\left[\D_{T_k}^-\right],A_k\right],
\end{equation}\abel{berg}
where we note that $\psi\left(\om,H_{\sigma_k}\right)$ is well-defined on $A_k$ since $\si_k<\infty$ on $A_k$.
Also observe that on $A_k$, $\psi\left(\om,H_{\sigma_k}\right)$ differs only at finitely many sites from $\om$ and 
therefore satisfies (\ref{pine}) since $\om$ does so. Hence we may use Lemma \ref{expect} 
and $\D_{T_k}^-\leq \D_{T_k}$ to conclude from (\ref{ber}) and (\ref{berg}) that
$E_{0,\om}[\Delta_k^-]\leq kP_{0,\om}[A_k]$. Therefore, due to (\ref{tuch}),
\[0=\lim_{K\to\infty}\frac{1}{K}\sum_{k=1}^KE_{0,\om}[\Delta_k^-]=\lim_{K\to\infty}\frac{E_{0,\om}\left[\D_{T_K}^-
\right]}{K}.\]
The claim now follows from
$E_{0,\om}\left[\D_{T_K}\right]=K$, see Lemma \ref{expect}, and $\D_{T_K}=\D_{T_K}^++\D_{T_K}^-$.
\end{proof}

Lemma \ref{eins}  and $\D_{T_K}^+\leq \sum_{x=0}^{K-1}\delta^{x}$ 
imply the following sufficient criterion for recurrence.
\begin{cor}\label{recdet}
0 is $\om$-recurrent if
\[
\liminf_{K\to\infty}\frac{1}{K}\sum_{x=0}^{K-1}\delta^{x}(\om)<1
\]
\end{cor}\abel{recdet}

We conclude this section by showing that the probability of never returning
to the starting point is positive whenever the starting point is $\om$-transient. In Section
\ref{leq2} we shall explicitly compute this probability in some cases.
\begin{lemma}\label{posse}
If 0 is $\om$-transient  then
 $P_{0,\om}[\forall n>0:\ X_n> 0]>0$.
\end{lemma}\abel{posse}
\begin{proof}
Since 0 is $\om$-transient  and since all positive integers are $P_{0,\om}$-a.s.\ 
eventually hit by the walk, we have $P_{0,\om}$-a.s.\ $X_n>0$ for $n$ large.
Now we distinguish two cases.

If 1 is $\om$-recurrent then it follows 
from the divergence part of the Borel Cantelli lemma that $\sum_{i}(1-\om(1,i))<\infty.$
Since 0 is $\om$-transient this  implies that
\begin{eqnarray*}
P_{0,\om}[\forall n>0:\ X_n> 0]&=& 
\om(0,1) P_{1,\om}\left[\forall n\ (X_n=1\Rightarrow X_{n+1}=2)\right]\\
&\geq& 
\om(0,1)\prod_{i\geq 0}\om(1,i)\ >\ 0
\end{eqnarray*}
as required.

If 1 is $\om$-transient  then there is some $K>0$ and a nearest neighbor path
$(x_n)_{n=0}^K$ of integers with $x_0=0$ and $x_K=2$ such that
\[
0<P_{0,\om}[H_K=(x_n)_{n=0}^K,\ \forall n\geq K:\ X_n\geq 2]=UV,
\]
where
\[ U:=P_{0,\om}[H_K=(x_n)_{n=0}^K]\quad\mbox{and }\quad V:=P_{2,\psi\left(\om,(x_n)_{n\leq K}\right)}[
\forall n\geq 0:\ X_n\geq 2].
\]
Hence $U>0$ and $V>0$.
Now we are doing some surgery on $(x_n)_{n=0}^K$ by cutting out the excursions from 1 toward 0.
To this end, we let
$n_1,\ldots,n_{k-1}$ be the enumeration of the times $1\leq n< K$ for which 
$x_n\geq 2$ or $x_{n+1}\geq 2$ and set $n_k:=K$. Then  $y_0:=0$ and $y_i:=x_{n_i}$ for $1\leq i\leq k$ 
defines 
  a nearest neighbor path $(y_n)_{n=0}^k$, which
starts at 0, ends at 2, and is strictly positive in between.  
Therefore,
\begin{eqnarray*}\lefteqn{
P_{0,\om}[\forall n>0:\ X_n> 0]}\\
&\geq& 
P_{0,\om}[H_k=(y_n)_{n=0}^k,\ \forall n\geq k:\ X_n\geq 2]=uv,
\end{eqnarray*}
where
\[u:=P_{0,\om}[H_k=(y_n)_{n=0}^k]\quad\mbox{ and }\quad
v:=P_{2,\psi\left(\om,(y_n)_{n\leq k}\right)}[\forall n\geq 0:\ X_n\geq 2].
\]
Consequently, it suffices to show that $u>0$ and $v>0$.
By Lemma \ref{same}, $V=v$
because
 $\psi\left(\om,(x_n)_{n\leq K}\right)(z)=\psi\left(\om,(y_n)_{n \leq k}\right)(z)$
 for $z\geq 2$ since $(y_n)_{n\leq k}$ visits each number $\geq 2$
as often as $(x_n)_{n\leq K}$ does. Hence, 
$v>0$ since $V>0$.

As for $U$ and $u$, both are products of finitely many 
 factors of the form $\om(x,i)$ and $1-\om(x,i)$. 
We have to make sure that none of the factors involved in $u$ 
is 0. Since $\om(x,i)\geq 1/2$, only terms of the form $1-\om(x,i)$ are critical.
Having a factor $1-\om(x,i)=0$ in $u$, which is not present in  $U$,
means that the path $(x_n)_n$ jumps to $x+1$ after the $i$-th visit to $x$
whereas $(y_n)_n$ jumps to $x-1$ after the $i$-th visit to $x$.
Since there are no steps from $1$ to 0 in $(y_n)_n$, any such $x$ 
must be at least 2. However, for any $x\geq 2$ all the steps from $x$ to $x+1$ and
from $x$ to $x-1$ happen in the same order for $(x_n)_n$ as for $(y_n)_n$, thus giving rise 
to the same factors $\om(x,i)$ and $1-\om(x,i)$ in $U$ and $u$. Consequently,
$u>0$ since $U>0$.
\end{proof}

\section{Recurrence and Transience in Random Environments}\label{ratire}\abel{ratire}
For $\PP$  a probability measure on $\Omega_+$, equipped with its canonical $\sigma$-field, and for $x\in\Z$ we define 
the semi-direct product $P_x:=\PP\times P_{x,\om}$ on $\Omega_+\times \Omega$  by
$P_x[\cdot]:=\EE\left[P_{x,\om}[\cdot]\right]$. This is the
 so-called \textit{annealed} measure which we get after averaging
the \textit{quenched} measure $P_{x,\om}$ over $\PP$. 
Here the expectation operators for $\PP$ and $P_x$ are denoted by $\EE$ and $E_x$,
respectively.

Not much can be said about recurrence and transience for general $\PP$. 
A conclusive answer can be given if
 $(\om(x))_{x\geq 0}$ is stationary and ergodic under $\PP$ w.r.t.\ the shift on $\Z$.
 Stationarity of $(\om(x))_{x\geq 0}$ means that  
the distribution of  $f\left(\theta^x(\om)\right)$ under $\PP$ for $x\geq 0$ does not depend on $x$, 
 where $f:\Omega_+\to\Omega_+$ is defined by 
\[
(f(\om))(x):=\left\{\begin{array}{ll}
\om(x)&\mbox{if $x\geq 0$}\\
1/2 &\mbox{if $x< 0$}
\end{array}\right.
\]
and  $\theta^x:\Omega_+\to\Omega_+$ is the canonical shift of $\om$ 
to the left by $x\ (x\in\Z)$ steps as defined by $(\theta^x(\om))(z):=\om(z+x)$.

\begin{rem}\label{rem2}{\rm In the special case where $\om(x,i)$ is $\PP$-a.s.\ for all $x\in\Z$  constant in $i$ (but not necessarily
constant in $x$), we get a one-dimensional \textit{random walk in random environment} (RWRE) with a nonnegative
drift. The general model of RWRE for $d=1$, which allows positive and negative drifts, 
has been studied e.g.\  by Solomon \cite{so}, see also \cite{sz} and \cite{zei} for results and references.
For a unifying model which includes RWRE and ERW 
we would have to replace $\Omega_+$ by $\Omega_{\pm}:=\left([-1,1]^\N\right)^\Z$. Our methods
do not immediately work in this case.}\hfill $\Box$
\end{rem}
\begin{theorem}\label{sun}
If $(\om(x))_{x\geq 0}$ is stationary and ergodic under $\PP$ then
either every $x\geq 0$ is $\PP$-a.s.\ $\om$-recurrent or 
       every $x\geq 0$ is $\PP$-a.s.\ $\om$-transient.
\end{theorem}\abel{sun}
In the first case mentioned above, i.e.\ when every $x\geq 0$ is $\PP$-a.s.\ $\om$-recurrent, 
 we shall call $(X_n)_n$ \textit{recurrent},
in the second case $(X_n)_n$  is called \textit{transient}. 
\begin{proof}
For all $x\geq 0$ and all $\om\in\Omega_+$ by Lemma \ref{fliege},  
\begin{equation}\label{kille}
P_{0,\om}[R_0]\leq P_{0,\om}[R_x]\ =
P_{-x,\theta^x(\om)}[R_0]  \stackrel{(\ref{wolfgang})}{=}  P_{0,f\left(\theta^x(\om)\right)}[R_0]. 
 \end{equation}\abel{kille}
Consequently, taking $\EE$-expectations in (\ref{kille}) and using stationarity yields
\[P_0[R_0]\ \leq\ \EE\left[P_{0,f\left(\theta^x(\om)\right)}[R_0]\right]\ =\ 
\EE\left[P_{0,f(\om)}[R_0]\right]\stackrel{(\ref{wolfgang})}{=}
\EE\left[P_{0,\om}[R_0]\right]\ =\ P_0[R_0].\]
Therefore, the inequality in (\ref{kille}) is in fact $\PP$-a.s.\ an equality. Hence,
$P_{0,f\left(\theta^x(\om)\right)}[R_0]$ does $\PP$-a.s.\ not depend on $x$. Moreover, the sequence
$P_{0,f\left(\theta^x(\om)\right)}[R_0]$\ $(x\geq 0)$ is ergodic because it is of the
form $g\left((\om(y))_{y\geq x}\right)\ (x\geq 0)$. Consequently, this sequence is $\PP$-a.s.\ 
equal to a deterministic constant, which is either 0 or 1 by Proposition \ref{rectrans}. 
 \end{proof}

The following lemma shows how the path inherits stationarity and/or
ergodicity from the environment.
\begin{lemma}\label{shift}
If $(\om(x))_{x\geq 0}$ is stationary (resp.\ ergodic) under $\PP$ then  
\[\xi:=(\xi_k)_{k\geq 0}:=\left(\left(\om(x+k)\right)_{x\geq 0},\left(X_{\tau_{k,m}}-k\right)_{m\geq 0}
\right)_{k\geq 0}\]
 is stationary (resp.\ ergodic)  under $P_0$.
In particular, for any measurable function $g$ on $([1/2,1]^{\N})^{\N_0}\times \Z^{\N_0}$ is the sequence
$(g(\xi_k))_{k\geq 0}$ stationary (resp.\ ergodic) under $P_0$ if  the sequence
$(\om(x))_{x\geq 0}$ is so under $\PP$.
\end{lemma}\abel{shift}
Here $\xi_k$ consists of the environment to the right of $k$ and of the part 
of the trajectory to the right of $k$.
\begin{proof}
To prove stationarity of  $\xi$  we shall show that for all measurable
subsets $B$ of the codomain of  $\xi$, $P_0\left[(\xi_k)_{k\geq K}\in B\right]$ is the same
for all $K\geq 0$.
For the proof of ergodicity we need to show that $P_0[A]\in\{0,1\}$ whenever there is a $B$ as above
such that 
\begin{equation}\label{AA}
A=\{(\xi_k)_{k\geq K}\in B\}\qquad\mbox{ for all $K\geq 0$.}
\end{equation}\abel{AA}
In both proofs the following identities will be used. For all $\om\in\Omega_+$, 
$K\geq 0$, and $B$ as above we have by the strong Markov property  $P_{0,\om}$-a.s.\
\begin{eqnarray}\label{glasses}
P_{0,\om}\left[(\xi_k)_{k\geq K}\in B\mid\F_{T_K}\right]
&=&
P_{K,\psi\left(\om,H_{T_K}\right)}\left[(\xi_k)_{k\geq K}\in B\right]\\ 
&=& \nonumber
P_{0,\theta^K\left(\psi\left(\om,H_{T_K}\right)\right)}\left[(\xi_k)_{k\geq 0}\in B\right].
\end{eqnarray}\abel{glasses}
Since $\theta^K\left(\psi\left(\om,H_{T_K}\right)\right)(x)$ and $f\left(\theta^K(\om)\right)(x)$ 
coincide $P_{0,\om}$-a.s.\
for $x\geq 0$ 
we can apply Lemma \ref{same} to see that (\ref{glasses}) equals
\begin{equation}\label{squirrel}
P_{0,f(\theta^K(\om))}\left[(\xi_k)_{k\geq 0}\in B\right]=:\eta_K(\om).
\end{equation}\abel{squirrel}
Hence taking $E_{0,\om}$-expectations in (\ref{glasses}) and (\ref{squirrel}) yields
\begin{equation}\label{ok}
P_{0,\om}\left[(\xi_k)_{k\geq K}\in B\right]=\eta_K(\om).
\end{equation}\abel{ok}
If we now take $\EE$-expectations on both side of (\ref{ok}) we get
\begin{eqnarray}\label{oaks}
P_0\left[(\xi_k)_{k\geq K}\in B\right]&=&\EE\left[P_{0,f(\theta^K(\om))}\left[(\xi_k)_{k\geq 0}\in B\right]\right]\\
&=&
\EE\left[P_{0,f(\om)}\left[(\xi_k)_{k\geq 0}\in B\right]\right],\nonumber
\end{eqnarray}\abel{oaks}
if $(\om(x))_{x\geq 0}$ is stationary 
under $\PP$. 
Hence in this case the left hand side of (\ref{oaks}) does not depend on $K$,
which  proves stationarity of $(\xi_x)_{x\geq 0}$.

For the proof of ergodicity of $\xi$ we assume (\ref{AA}).
Then $\eta_K$ does not depend on $K$ since  the left-hand side  of 
(\ref{ok})  does not. However, 
since $\eta_K$ is a function of the form  $g((\om(x))_{x\geq K})$, 
the process $(\eta_K)_{K\geq 0}$ is ergodic 
if  $(\om(x))_{x\geq 0}$ is so. Therefore, in this case, being independent of $K$, $\eta_K$ 
is  $\PP$-a.s.\ equal to a deterministic
constant $c$. Going back from (\ref{ok}) via (\ref{squirrel}) to (\ref{glasses}) we obtain that
$\PP$-a.s.\ 
$P_{0,\om}\left[A\mid\F_{T_K}\right]=c$.
However, given $\om$, $A\in\sigma\left(\bigcup_K \F_{T_K}\right)$. Therefore,
 by the martingale convergence theorem, $P_{0,\om}$-a.s.\ $P_{0,\om}\left[A\mid\F_{T_K}\right]
\to\won_A$
as $K\to\infty$. 
Hence, either $P_0$-a.s.\  $c=
P_{0,\om}\left[A\mid\F_{T_K}\right]=0$  or $P_0$-a.s.\
 $c=
P_{0,\om}\left[A\mid\F_{T_K}\right]=1$.
Integration w.r.t.\ $P_0$ gives $P_0[A]\in\{0,1\}$.
\end{proof}
The next result deals with the maximal cookie consumption per site.
\begin{lemma}\label{less}
If $(\om(x))_{x\geq 0}$ is stationary under $\PP$ then $E_0[\D_\infty^x]\leq 1$ for all $x\geq 0$. 
\end{lemma}\abel{less}
\begin{proof}
Let $x\geq 0$. For $\om$ given, the distribution of 
$\D_\infty^x$ depends only on the distribution of $(X_{\tau_{0,m}})_{m\geq 0}$.
Therefore,
 we may assume due to Lemma \ref{same} without loss of generality, that (\ref{pine}) is $\PP$-a.s.\ fulfilled. 
Let $0\leq k< K$. Then $P_0$-a.s.\
\[\D_{T_K}\ \geq\ \D_{T_K}^+=\sum_{y=0}^{K-1}\D_{T_K}^y\ \geq\ \sum_{y=0}^{K-1-k}\D_{T_K}^y\ \geq\ 
\sum_{y=0}^{K-1-k}\D_{T_{y+k}}^y
\]
since $T_{y+k}< T_{K}$ for $y<K-k$. 
Consequently, we obtain from Lemma \ref{expect},
\[K=E_0[\D_{T_K}]\ \geq\  \sum_{y=0}^{K-1-k}E_0\left[\D_{T_{y+k}}^y\right]\ =\ 
(K-k)E_0\left[\D_{T_{x+k}}^x\right]
\]
 due to stationarity of the sequences $(\D_{T_{y+k}}^y)_{y\geq 0}$\ $(k\geq 0)$, which we get from 
  Lemma \ref{shift} applied for all $k\geq 0$ to
\[g\left((\om(x))_{x\geq 0}, (x_m)_{m\geq 0}\right):=\sum_{i=1}^{\#\{m<T_k((x_n)_n)\mid x_m=0\}}
(2\om(0,i)-1).\]
 Therefore, $E_0[\D_{T_{x+k}}^x]\leq K/(K-k)$. Letting $K\to\infty$ gives
$E_0[\D_{T_{x+k}}^x]\leq 1$ for all $k\geq 0$. Monotone convergence as $k\to\infty$ then yields the claim.
\end{proof}

The second part of the following theorem  classifies recurrent and transient walks.
\begin{theorem}\label{main}
Assume that $(\om(x))_{x\geq 0}$ is stationary and ergodic under $\PP$.
Then 
\begin{equation}\label{equal}
E_0[\D_{\infty}^x]=\min\left\{1,\EE[\delta^0]\right\}\quad\mbox{ for all $x\geq 0$.}
\end{equation}
Moreover,
if 
\begin{equation}\label{witri}
\PP[\om(0)=(1,1/2,1/2,1/2,\ldots)]<1
\end{equation}\abel{witri} then
\begin{equation}\label{label}
\mbox{$(X_n)_{n\geq 0}$ is recurrent if and only if
$\EE\left[\delta^0\right]\leq 1.$}
\end{equation}\abel{label}
\end{theorem}\abel{main}\abel{equal}
Obviously, 
if (\ref{witri}) fails then 
$X_n=n$ $P_0$-a.s.\ for all $n$, which makes the walk transient although $\EE[\delta^0]=1$.
\begin{proof}
 Lemma \ref{shift} applied to
\[g\left((\om(x))_{x\geq 0}, (x_m)_{m\geq 0}\right):=\sum_{i=1}^{\#\{m\mid x_m=0\}}
(2\om(0,i)-1)\]
yields that $(\D_\infty^k)_{k\geq 0}$ is stationary.
Therefore,  we may assume for the
proof of (\ref{equal}) that  $x=0$.
Moreover, since $\won_{R_0}$ and $\D_{\infty}^0$  are functions of $(\om(x))_{x\geq 0}$ and $(X_{\tau_0,m})_{m\geq 0}$
 we may 
assume thanks to Lemma \ref{same}  without loss of generality
that the assumption (\ref{pine})
is satisfied for all $\om\in\Omega_+$.

Due to Theorem \ref{sun}, $(X_n)_n$ is either recurrent or transient.
If it is recurrent
then the walker will eat $P_0$-a.s.\ all the cookies at $0$, which results 
in $\D_\infty^0= \delta^0$, thus showing
 $E_0[\delta^0]=E_0[\D_\infty^0]$.  
Lemma \ref{less} then yields $\EE\left[\delta^0\right]\leq 1$ and
(\ref{equal}).

Now we assume that the walk  is transient. Then by 
stationarity of $(\D_\infty^k)_{k\geq 0}$,  see above,
\[
E_0[\D_\infty^0]=\frac{1}{K}\sum_{k=0}^{K-1}E_0[\D_\infty^k]\ =\
\EE\left[\frac{1}{K}E_{0,\om}\left[\sum_{k=0}^{K-1}\D_\infty^k\right]\right]
\ \geq\ \EE\left[\frac{E_{0,\om}\left[\D_{T_K}^+\right]}{K}\right].
\]
Since $E_{0,\om}\left[\D_{T_K}^+\right]/K\leq E_{0,\om}\left[\D_{T_K}\right]/K=1$,  see
 Lemma \ref{expect}, we get by  dominated convergence and Lemma \ref{eins} that
$E_0[\D_\infty^0]\geq 1$. Since the opposite inequality holds due to  Lemma \ref{less} we conclude
\begin{equation}\label{erlangen}
E_0[\D_\infty^0]= 1.
\end{equation}\abel{erlangen}
Now consider the event 
\[ S:=\left\{\sum_{i\geq 2}(2\om(0,i)-1)>0\right\}\]
that 0 has not all  its drift stored in its first cookie.
We claim $\PP[S]>0$. Indeed, otherwise $1>\EE[\delta^0]$
since we 
 excluded the degenerate case in which the first cookie has $\PP$-a.s.\ parameter 1.  
 However, because of $\delta^0\geq \D_\infty^0$ this would contradict (\ref{erlangen}).

By Lemma \ref{posse}, 
$\PP$-a.s.\ $P_{0,\om}[\forall n>0:\ X_n> 0]>0$.
Therefore, since  $\PP[S]>0$ we have 
\begin{eqnarray*}
0&<&\EE\left[P_{0,\om}[\forall n>0:\ X_n> 0], S\right]=P_0\left[\{\D_\infty^0=2\om(0,1)-1\}\cap S\right]\\
&\leq& P_0[\D_\infty^0<\delta^0].
\end{eqnarray*}
Since $\D_\infty^0\leq\delta^0$ this implies  $E_0[\D_\infty^0]<E_0[\delta^0]$. 
From this we obtain
by (\ref{erlangen}) that $1<\EE[\delta^0]$ and (\ref{equal}) as required.
\end{proof}

\section{Eating left-overs}\label{elo}\abel{elo}
Assume that $(\om(x))_{x\geq 0}$ is stationary and ergodic. Then by Theorem
\ref{sun} $(X_n)_{n\geq 0}$ is either recurrent or transient and Theorem
\ref{main} tells us which is the case. 

Let us  assume that   $(X_n)_{n\geq 0}$ is transient. Then the walk will visit each site $x\in\Z$
$P_0$-a.s.\ only a finite number of times. Hence $\om_2:=\psi\left(\om,
(X_n)_{n<\infty}\right)$, with a straightforward extension of definition (\ref{psi})
to infinite sequences, is $P_0$-a.s.\ well defined and consists of the
cookies left over by the random walk.  
So we may start a second ERW $(X_n^{(2)})_{n\geq 0}$ in the environment
$\om_2$. Since $(\om(x))_{x\geq 0}$ was stationary and ergodic, so is
$(\om_2(x))_{x\geq 0}$ due to Lemma \ref{shift} applied to
\[g\left(\om(x))_{x\geq 0}, (x_m)_{m\geq 0}\right):=\left(\om\left(0,i+\#\{m\geq 0\mid x_m=0\}\right)
\right)_{i\geq 1}.
\]
Consequently, also $(X_n^{(2)})_{n\geq 0}$ is either recurrent or transient.
Moreover, due to (\ref{equal}) and (\ref{label})  the first random walk has reduced the expected total 
drift stored in the cookies at any site $x\geq 0$ by 1. 
If the total drift stored in the cookies, which were  left over by the first walk,
 is less than 1 then  
$(X_n^{(2)})_{n\geq 0}$ will be recurrent due to Theorem \ref{main}. If it is larger than
1 then it will be transient and will  leave behind another stationary and ergodic
 environment 
$\om_3:=\psi(\om_2,
(X_n^{(2)})_{n<\infty})$, in which we can start a third ERW
$(X_n^{(3)})_{n\geq 0}$. 

This can be iterated.
E.g.\ if $\EE[\delta^0]$ is finite but not an integer (to avoid exceptions related to the one 
ruled out in (\ref{witri})) then the first $\lfloor \EE[\delta^0]\rfloor$ ERWs
will almost surely be transient and the next one will almost surely be recurrent and
 will eventually eat all
the cookies on $\N_0$.
\section{Strong law of large numbers}\label{sloln}\abel{sloln}

\begin{theorem}\label{lln}
If $(\om(x))_{x\geq 0}$ is stationary and ergodic under $\PP$ then $P_0$-a.s.\
\[\lim_{n\to\infty}\frac{X_n}{n}=v:=\frac{1}{u}\geq 0,\quad\mbox{where}\quad
u:=\sum_{j\geq 1}P_0\left[T_{j+1}-T_j\geq j\right]\in[1,\infty].\]
\end{theorem}\abel{lln}
Roughly speaking, $u$ is the expected time it takes a walker who has just arrived at $\infty$
to reach level $\infty+1$. One could phrase the proof of Theorem \ref{lln} in terms
of a limiting distribution of the environment viewed from the particle. The following proof
is a bit more elementary.
\begin{proof}
We shall show that $P_0$-a.s.\ $T_k/k\to u$  
as $k\to\infty$.
It is  standard (see e.g.\ \cite[Lemma 2.1.17]{zei})
that this implies that $X_n/n$ converges $P_0$-a.s.\ to $1/u$.
Observe that $P_0$-a.s.
\begin{equation}\label{stat}
T_k=\sum_{i=0}^{k-1}T_{i+1}-T_i.
\end{equation}\abel{stat}
Consequently,
\begin{equation}\label{kloss}
\liminf_{k\to\infty}\frac{T_k}{k}\geq\sup_{t\geq 0}\liminf_{k\to\infty}\frac{1}{k}\sum_{i=t}^{k-1}
\left((T_{i+1}-T_i)\wedge t\right).
\end{equation}\abel{kloss}
Applying  Lemma \ref{shift} for all $t\geq 0$ to
\[g\left((\om(x))_{x\geq 0},(x_m)_{m\geq 0}\right):=(T_{t+1}-T_t)\left((x_m)_{m\geq 0}\right)\wedge t\]
yields
that $\left((T_{i+1}-T_i)\wedge t\right)_{i\geq t}$ is stationary
and ergodic for all $t\geq 0$.
Therefore, by the ergodic theorem, the right-hand side of (\ref{kloss}) is equal to
\begin{eqnarray*}
&&\sup_{t\geq 0}E_0\left[(T_{t+1}-T_t)\wedge t\right]
\ =\ \nonumber
\sup_{t\geq 0}\sum_{j= 1}^tP_0\left[(T_{t+1}-T_t)\wedge t\geq j\right]\nonumber
\\
&=&
\sup_{t\geq 0}\sum_{j= 1}^tP_0\left[(T_{t+1}-T_t)\wedge j\geq j\right]\ =\
\sup_{t\geq 0}\sum_{j= 1}^tP_0\left[T_{j+1}-T_j\geq j\right]\ =\ u,
\end{eqnarray*}
where we used in the second to last inequality stationarity of 
$\left((T_{i+1}-T_i)\wedge j\right)_{i\geq j}$ for all $j$.
On the other hand, (\ref{stat}) implies
\begin{eqnarray}
\limsup_{k\to\infty}\frac{T_k}{k}&=&
\limsup_{k\to\infty}\frac{1}{k}\sum_{i=0}^{k-1}\sum_{j\geq 1}\won_{T_{i+1}-T_i\geq j}
\nonumber\\
&\leq &\sum_{j\geq 1}\limsup_{k\to\infty}\frac{1}{k}\sum_{i=0}^{k-1}\won_{T_{i+1}-T_i\geq j}.
\label{quak}
\end{eqnarray}\abel{quak}
Due to Lemma \ref{shift} applied to
\[ g\left((\om(x))_{x\geq 0},(x_m)_{m\geq 0}\right):=\won\left\{(T_{j+1}-T_j)\left((x_m)_{m\geq 0}\right)
\geq j\right\}\]
the sequences $\left(\won_{T_{i+1}-T_i\geq j}\right)_{i\geq j}$,\ $j\geq 1,$
are stationary and ergodic. Consequently, by the ergodic theorem, 
the right-hand  side of (\ref{quak}) is $P_0$-a.s.\ equal to
 $u$, too.
\end{proof}

\begin{rem}\label{rem3}
{\rm The one-dimensional model under consideration can be extended to  higher
dimensions $d\geq 1$ by letting $\om(x,i)\in[0,1]^{2d}$\ $(x\in\Z^d, i\geq 1)$ be a 
vector of transition probabilities to the $2d$ neighbors of $x$ in $\Z^d$. In the case where
$\om(x),\ x\in\Z^d,$ are i.i.d.\ under $\PP$, a straightforward adaptation of a renewal
structure technique introduced by Sznitman and Zerner for random walks
in random environments (RWRE) gives the $P_0$-a.s.\ convergence of $(X_n\cdot\ell)/n$
towards a deterministic limit on the event $\{\lim_{n\to\infty} X_n\cdot\ell=+\infty\}$,
where $\ell$ is any direction in $\R^d$, see  \cite{szze} and \cite[Theorem 3.2.2]{zei}.
Here we assume that none of the transition probabilities is equal to 0. In this case
Lemma \ref{posse} and its higher dimensional analogue (see \cite[(1.16)]{szze}) are easy to obtain.
}\hfill $\Box$
\end{rem}

\section{Monotonicity}\label{m}
\begin{figure}\label{cup}
\epsfig{file=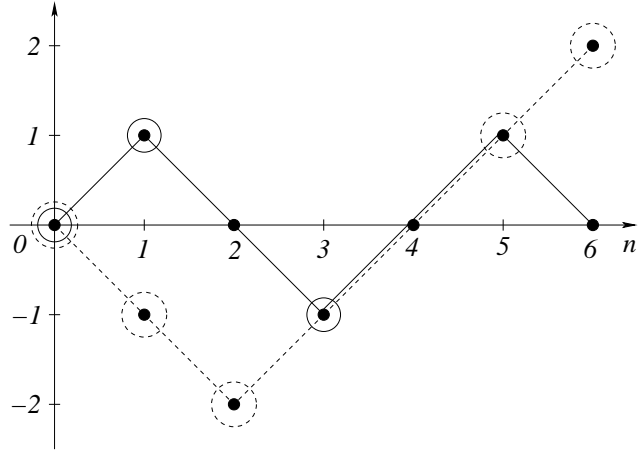,height=6cm}
\caption{How  $(X_n^{(1)})_n$ (dashed), which eats weak cookies, might overtake 
 $(X_n^{(2)})_{n}$ (solid), which  eats strong cookies, 
see Example \ref{standard}. 
Sites which are visited for the first time are encircled accordingly.
}
\end{figure}

\abel{m}
Monotonicity results are often difficult to obtain for processes in random media
since standard coupling techniques, similar to the one used in the proof of Lemma \ref{simply} 
and Example \ref{standard}, see below, tend to fail.
The following result shows that, roughly speaking, starting further to the right,
helps to reach a goal located to the right sooner.
\begin{lemma}{\rm (Monotonicity w.r.t.\ initial point)}\label{monoinitial}
Let $\om\in\Omega_+$, 
$-\infty\leq x\leq y_1\leq y_2\leq z\leq \infty$,\ $y_1, y_2\in\Z,$
and $t\in[0,\infty]$. Then 
\begin{equation}\label{grille}
P_{y_1,\om}[T_z\leq T_x\wedge t]\leq P_{y_2,\om}[T_z\leq T_x\wedge t].
\end{equation}\abel{grille}
\end{lemma}\abel{monoinitial}
Here we define $T_{\infty}=T_{-\infty}=\infty$.
\begin{proof}
By continuity it is enough to show the claim for $z<\infty$.
Moreover, by induction it suffices to show the statement for  $y_2=y_1+1$. So assume $y_2=y_1+1$.
For $y_1<z$ 
denote by $\Pi_{y_1}^z$ the set of all finite nearest-neighbor paths $\pi=(x_n)_{n\leq m},$\ $m> 0,$ 
which  start at $x_0=y_1$, end at $x_m=z$ and do not hit $z$ in between. 
Any such path $\pi$ can be uniquely written as the concatenation
$\pi=(B_1, A_1, B_2, A_2,\ldots,B_{j(\pi)}, A_{j(\pi)})$ for some $j(\pi)\geq 1$,
where $A_i$ and $B_i$ $(i\leq j(\pi))$
are nonempty nearest-neighbor paths such that the $A_i$'s contain only points $> y_1$ (``{\bf A}bove $y_1$")
and the  $B_i$'s contain only points $\leq y_1$ (``{\bf B}elow $y_1$").
Then the function $\Phi: \Pi^z_{y_1}\to\Pi^z_{y_1+1}$ defined by
\[\Phi(B_1, A_1, B_2, A_2,\ldots,B_{j}, A_{j}):=(A_1, B_1, A_2, B_2,\ldots,
A_{j-1},B_{j-1}, A_{j}),\]
  see Figure \ref{surgery}, is well-defined and surjective.
\begin{figure}\label{surgery}
\epsfig{file=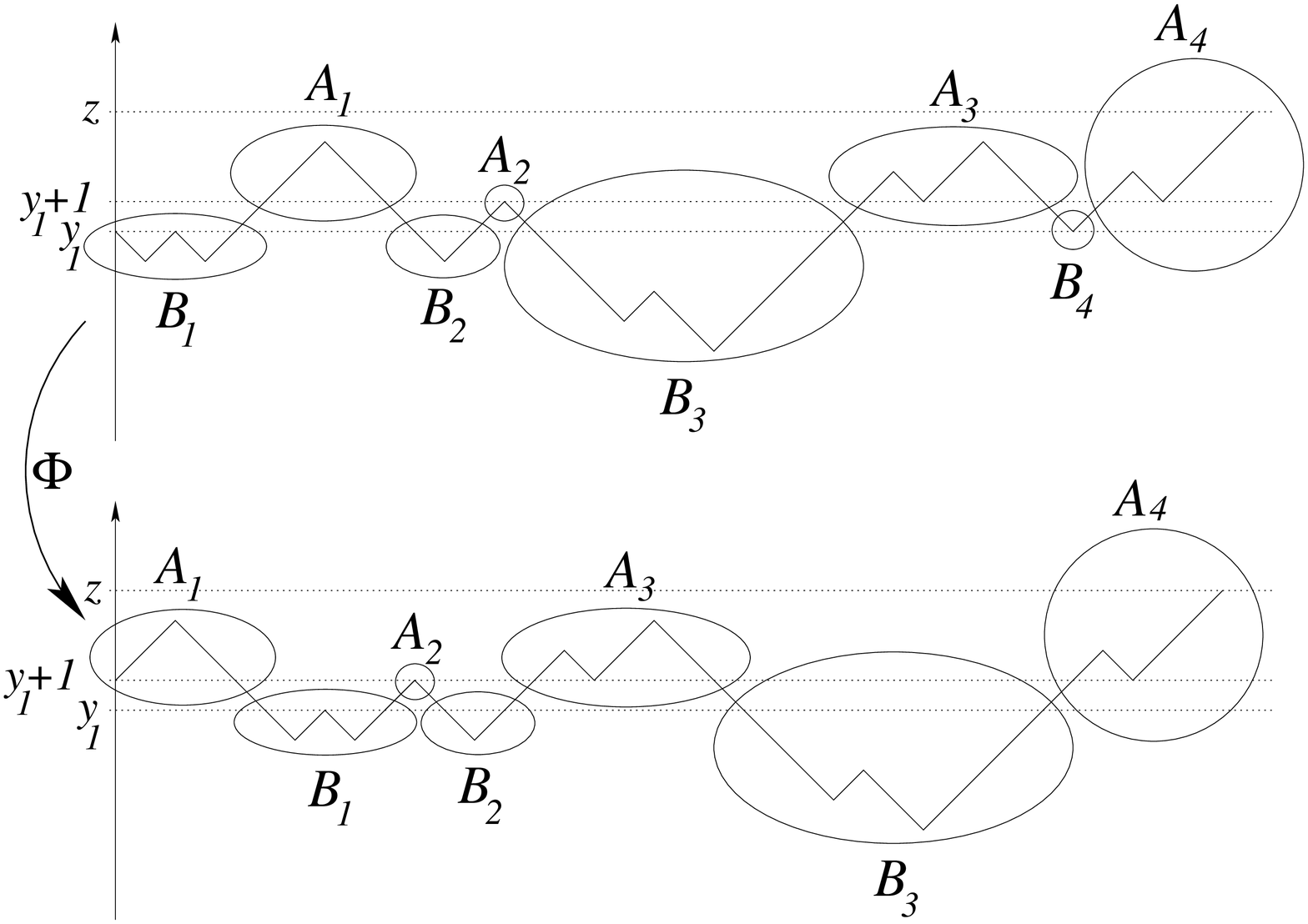,height=6cm,width=11cm}
\caption{Transforming paths from $y_1$ to $z$ into paths from $y_1+1$ to $z$. 
}
\end{figure}\abel{surgery}
This function cuts out the last excursion $B_j$ from $y_1$ downward  but otherwise
only rearranges pieces of $\pi$ without changing the relative order
in which the points above $y_1$ are visited nor the  relative order
in which the points below $y_1$ are visited. 
Therefore, for any $\pi\in\Pi^z_{y_1}$,
\[P_{y_1,\om}\left[H_{T_z}=\pi\right]=
P_{y_1+1,\om}\left[H_{T_z}=\Phi(\pi)\right]\cdot 
P_{y_1,\om'}\left[H_{T_{y_1+1}}=(B_{j(\pi)},y_1+1)\right],
\]
where
$\om':=\psi(\om,(B_1, A_1,\ldots,B_{j(\pi)-1}, A_{j(\pi)-1},y_1))$
is the environment  $\pi$ faces before it starts its excursion $B_{j(\pi)}$.
By summing over all possible excursions $B_{j(\pi)}$ we get for all $\pi\in\Pi^z_{y_1}$,
\[P_{y_1,\om}\left[H_{T_z}\in \Phi^{-1}(\{\Phi(\pi)\})\right]=
P_{y_1+1,\om}\left[H_{T_z}=\Phi(\pi)\right].
\]
Since $\Phi$ is surjective this means that  for all $\pi\in\Pi^z_{y_1+1}$,
\begin{equation}\label{aufpassen}
P_{y_1,\om}\left[H_{T_z}\in \Phi^{-1}(\{\pi\})\right]=
P_{y_1+1,\om}\left[H_{T_z}=\pi\right].
\end{equation}\abel{aufpassen}
Now denote by $\Pi_{y_1}^z(x,t)$ the set of paths $\pi$ in $\Pi_{y_1}^z$ which do not visit $x$ before
$z$ and which make at most
$t$ steps.
Then the right-hand side of (\ref{grille}) can be written 
as
\begin{equation}P_{y_1+1,\om}\left[H_{T_z}\in\Pi_{y_1+1}^z(x,t)\right]
=P_{y_1,\om}\left[H_{T_z}\in\Phi^{-1}\left(\Pi_{y_1+1}
^z(x,t)\right)\right]\label{pd}
\end{equation}\abel{pd}
due to (\ref{aufpassen}).
Cutting out an excursion does not make a path longer 
nor does it make a path visit $x$ if the path  did not do so before. Therefore,
\[\Phi^{-1}\left(\Pi_{y_1+1}
^z(x,t)\right)\supseteq \Phi^{-1}\left(\Phi\left(\Pi_{y_1}^z(x,t)\right)\right)\supseteq 
\Pi_{y_1}^z(x,t).\]
Consequently, the right-hand side
of (\ref{pd}) is greater than or equal to the left-hand side of (\ref{grille}).
\end{proof}
Roughly speaking, the following result states that increasing the strength of some cookies
does not slow down the walk. Here
we denote by  $\leq$ the canonical partial order on $\Omega_+$, i.e.\ $\om_1\leq\om_2$
 if and only if $\om_1(x,i)\leq \om_2(x,i)$ for all $x\in\Z$ and all $i\geq 1$.
\begin{lemma}{\rm (Monotonicity w.r.t.\ environment)}\label{monoenv}
Let $\om_1,\om_2\in\Omega_+$ with $\om_1\leq \om_2$ and $-\infty\leq x\leq y\leq z\leq +\infty$,\ $y\in\Z$,
and $t\in\N\cup\{\infty\}$. Then
\begin{equation}\label{ziept}
P_{y,\om_1}[T_z\leq T_x\wedge t]\leq P_{y,\om_2}[T_z\leq T_x\wedge t].
\end{equation}\abel{ziept}
\end{lemma}\abel{monoenv}
Intuitively, this result seems to be clear. However, 
the following example shows that the 
naive coupling approach to prove monotonicity  w.r.t.\ the environment
fails.
\begin{ex}\label{standard}
{\rm Let $\om_1,\om_2\in\Omega_+$ with $\om_j(x,1)=p_j$ and $\om_j(x,i)=1/2$ for $x\in\Z, i\geq 2$ and $j=1,2$,
where $1/2<p_1<p_2<1$. Thus $\om_1\leq \om_2$. There does not seem to be a simple way
to couple like in the proof of Lemma \ref{simply} two ERWs $(X_n^{(1)})_n$ and $(X_n^{(2)})_n$  in the environments 
$\om_1$ and $\om_2$, respectively, such that
$X_n^{(1)}\leq X_n^{(2)}$ for all $n\geq 0$ almost surely,  as we shall show now:

Again, let $U_n,\ n\geq 0,$ be a sequence of independent random variables
uniformly distributed on $[0,1]$.
If  the walk $(X_n^{(j)})_n$\ ($j=1,2$) visits at time $n$ a site for the first time then it moves to the
right in the next step iff $U_n<p_j$. If it has visited the site it is currently at at least once before 
than it moves to the right iff $U_n<1/2$.
Clearly, this defines two ERWs in the environments $\om_1$ and $\om_2$.
However,  $X_6^{(1)}>X_6^{(2)}$  on the event
\[\left\{(U_n)_{n=0}^5\in\ ]p_1,p_2[\ \times\ ]p_2,1]\ \times\ ]1/2,p_1[ \ \times\ [0,1/2[^2\ \times\ ]1/2,p_1[
\right\},
\]
which has positive probability, see Figure \ref{cup}. \
}\hfill $\Box$
\end{ex}\abel{standard}
\abel{cup}

\begin{proof}[Proof of Lemma \ref{monoenv}]
We assume that $t$ is finite. The case of $t=\infty$ follows then by continuity
 from the finite case.\\  By time $t$ the walker can eat only cookies which
are among the first $t$ cookies  at sites which are within distance $t$
from the starting point $y$.
Hence we may assume that $\om_1$ and $\om_2$ differ in the strength of
 only a finite number $r$ of cookies. By induction it suffices to consider the case $r=1$.
So let us assume that $\om_1(u,i)= \om_2(u,i)$ for all $(u,i)\in(\Z\times\N)\backslash
\{(v,j)\}$, where $(v,j)$ denotes location and number of the only cookie which might
be stronger in $\om_2$ than in $\om_1$. We shall refer to this cookie as to the crucial cookie.
 Since for the
event $\{T_z<T_x\wedge t\}$ only cookies between $x$ and $z$ matter we
may additionally assume $x<v<z$. 
Denote by $S$ the time of the $j$-th visit to $v$.  This is the time at which the walk
reaches the crucial cookie.
Then for any $i=1,2$,
\begin{equation}\label{ischias}
P_{y,\om_i}[T_z<T_x\wedge t]=P_{y,\om_i}[T_z<T_x\wedge t\wedge S]+P_{y,\om_i}[S<T_z<T_x\wedge t].
\end{equation}\abel{ischias}
Note that the first term on the right-hand side of (\ref{ischias}) does not depend on $i$.
It therefore suffices to show that the second term is non-decreasing in $i$.
 By the  Markov property the second term equals 
\begin{equation}\label{affe}
\sum_{s\geq 0}^{t-2}E_{y,\om_i}\left[P_{v,\psi\left(\om_i,H_{s}\right)}[
T_z<T_x\wedge (t-s)],\ 
s=S<T_z\wedge T_x\right].
\end{equation}\abel{affe}
Decomposition after the next  step, in which the strength of the crucial cookie comes into play,
yields $P_{y,\om_i}$-a.s.\ $(i=1,2)$
\begin{eqnarray}
\label{piep}\lefteqn{P_{v,\psi\left(\om_i,H_{s}\right)}[
T_z<T_x\wedge (t-s)]}\\
&=&(1-\om_i(v,j))P_{v-1,\psi\left(\om_i,H_{s+1}\right)}[
T_z<T_x\wedge (t-s-1)]\nonumber\\
&&+\ \om_i(v,j)P_{v+1,\psi\left(\om_i,H_{s+1}\right)}[
T_z<T_x\wedge (t-s-1)].\nonumber
\end{eqnarray}\abel{piep}
However, on $\{S=s\}$,
\begin{equation}\label{haiti}
\psi(\om_1,H_{s+1})=\psi(\om_2,H_{s+1})
\end{equation}\abel{haiti}
$P_{y,\om_1}$-a.s.\ because at  time $S$ the walker reaches and then eats the crucial 
cookie, thus erasing the only difference between
the two  environments. Moreover, we may apply Lemma \ref{monoinitial} with $y_1=v-1$ and $y_2=v+1$
to deduce that $P_{y,\om_1}$-a.s.
\[P_{v-1,\psi\left(\om_1,H_{s+1}\right)}[
T_z<T_x\wedge (t-s-1)]\leq P_{v+1,\psi\left(\om_1,H_{s+1}\right)}[
T_z<T_x\wedge (t-s-1)].
\]
Combined with (\ref{haiti}) and $\om_1(v,j)\leq \om_2(v,j)$ this implies that (\ref{piep}) for $i=1$ is $P_{y,\om_1}$-a.s.\
less than or equal to (\ref{piep}) for $i=2$. Consequently, (\ref{affe})  for $i=1$ is 
less than or equal to
\begin{equation}\label{tutu}
\sum_{s\geq 0}^{t-2}E_{y,\om_1}\left[P_{v,\psi\left(\om_2,H_{s}\right)}[
T_z<T_x\wedge (t-s)],\ 
s=S<T_z\wedge T_x\right].
\end{equation}\abel{tutu}
 However, 
the distribution of the above integrand under $E_{y,\om_1}$ does not depend on
$\om_1(v,j)$ any more. Therefore, we can replace $E_{y,\om_1}$ in expression (\ref{tutu}) by $E_{y,\om_2}$
without changing its value, thus getting (\ref{affe}) with $i=2$.
\end{proof}

The following corollaries show that the probability of return to the origin and the speed
of the walk are monotone increasing in $\om$.
\begin{theorem}\label{return}
The probability $P_{0,\om}[\forall n>0\ X_n> 0]$
never to return to the initial point is monotone increasing in $\om$.
\end{theorem}\abel{return}
\begin{proof}
By the simple Markov property for all $\om\in\Omega_+$,
\begin{eqnarray}\nonumber
P_{0,\om}[\forall n>0\ X_n> 0]&=&\om(0,1)P_{1,\psi(\om,(0,1))}[\forall n>0\  X_n>0]\\
&=&\om(0,1)P_{1,\om}[\forall n>0\  X_n>0],\label{simple}
\end{eqnarray}\abel{simple}
which is monotone increasing in $\om$ due to Lemma \ref{monoenv} applied to $x=0, y=1$ and $z=t=\infty$.
\end{proof}
\begin{theorem}\label{velo}
Let $\overline{\PP}$ be a probability measure on $\Omega_+^2$ such that 
\[\overline{\PP}\left[\{(\om_1,\om_2)\in \Omega_+^2\mid \om_1\leq \om_2\}\right]=1\] and
such that $(\om_i(x))_{x\geq 0}$ is stationary and ergodic for $i=1,2$.
Then $v_1\leq v_2$ if we denote by $v_i$\  $(i=1,2)$ the $\overline{\PP}\times P_{0,\om_i}$-a.s.\ constant limit of
$X_n/n$ as $n\to\infty$.
\end{theorem}\abel{velo}
\begin{proof}
By dominated convergence for $i=1,2$,
\begin{eqnarray*}
v_i&=&\lim_{k\to\infty}\overline{\EE}\left[E_{0,\om_i}\left[X_{T_k}/T_k
\right]\right]\ =\ 
\lim_{k\to\infty}\overline{\EE}\left[E_{0,\om_i}\left[k/T_k
\right]\right]\\
&=&\lim_{k\to\infty}\overline{\EE}\left[\int_0^1P_{0,\om_i}\left[T_k\leq k/t\right]\ dt\right].
\end{eqnarray*}
The statement now follows from Lemma \ref{monoenv} with $x=-\infty,\ y=0$, and $z=k$.  
\end{proof}

\begin{problem}
{\rm  Are the return probability and the velocity in an appropriate sense continuous in $\om$?
}
\end{problem}

\section{No excitement after the second visit}\label{leq2}\abel{leq2}
In some cases 
when the random walk behaves on the third and any later visit to a site like a simple symmetric
random walk
one can determine the probability that the walk will never  return to its starting point 
and can show that the  walk has zero speed.
\begin{theorem}\label{hp}
If $(\om(x))_{x\geq 0}$ is i.i.d.\ under $\PP$ such that
$\PP$-a.s.\ $\om(0,i)=1/2$ for all $i\geq 3$  and $\PP[\om(0,2)=1/2]<1$ then 
\begin{equation}\label{nani}
P_0[\forall n>0\ X_n>0]=\frac{\EE[\om(0,1)](\EE[\delta^0]-1)_+}{\EE[(2\om(0,2)-1)\om(0,1)]}.
\end{equation}
\end{theorem}\abel{hp}\abel{nani}
\begin{proof}
Consider $\delta^0-\D_\infty^0$, the total drift stored in the cookies at 0 which will never
be  eaten by
the random walk. On the one hand, by (\ref{equal})
\[
E_0[\delta^0-\D_\infty^0]=(\EE[\delta^0]-1)_+.
\]
On the other hand, since the first cookie at 0 is eaten $P_0$-a.s.\ right at the beginning of the walk and
since only the first two cookies at 0 contribute to $\delta^0$
we have $P_0$-a.s.\
\[
  \delta^0-\D_\infty^0=(2\om(0,2)-1)\won\{\forall n>0\  X_n>0\}.
\]
Combining these two facts we get
\begin{equation}\label{c3}
(\EE[\delta^0]-1)_+=\EE\left[(2\om(0,2)-1)P_{0,\om}[\forall n>0\  X_n>0]\right].
\end{equation}\abel{c3}
Recall (\ref{simple}) and note that
 $P_{1,\om}[\forall n>0\  X_n>0]$ is a function of $(\om(x))_{x\geq 1}$. Therefore,
it is independent of $\om(0)$ under $\PP$ by assumption. This has two consequences: Firstly,
taking $\EE$-expectations in (\ref{simple}) yields
\begin{equation}\label{money}
P_0[\forall n>0\  X_n>0]=\EE[\om(0,1)]P_1[\forall n>0\  X_n>0].
\end{equation}\abel{money}
Secondly, substituting (\ref{simple}) into (\ref{c3}) gives
\[
(\EE[\delta^0]-1)_+=\EE\left[(2\om(0,2)-1)\om(0,1)\right] P_{1}[\forall n>0\  X_n>0].
\]
Combined with 
 (\ref{money}) this proves the claim.
\end{proof}
\begin{theorem}\label{v0}
Let $(\om(x))_{x\geq 0}$ be stationary and ergodic with $\om(0,i)=1/2$ $\PP$-a.s.\
for all $i\geq 3$ and $\PP[\om(0,1)<1, \om(1,1)<1]>0$. Then
\[
\lim_{n\to\infty}\frac{X_n}{n}=0\quad\mbox{$P_0$-a.s..}
\] 
\end{theorem}\abel{v0}
\begin{proof}
By assumption, we can fix  $\eps>0$ such that  $\PP[A_{j}]=:\al$ is strictly positive
and independent of $j$, where
\[A_{j}:=\{\om(j-1,1)<1-\eps,\ \om(j,1)<1-\eps\}\quad (j\geq 1).\]
Due to  Theorem \ref{lln} we need to show that $u=\infty$.
To simplify
calculations we will do a worst case analysis by maximizing the strength
of selected cookies as follows.
Define $\bar{\om}_j\in\Omega_+$ for $j\in\Z$ by  
\[\bar{\om}_j(x):=\left\{\begin{array}{llll}
(1,&1/2,&1/2,1/2,\ldots)&\mbox{if $x<j-1$},\\
(1/2,&1/2,&1/2,1/2,\ldots)&\mbox{if $x=j-1$},\\
(1-\eps,&1,&1/2,1/2,\ldots)&\mbox{if $x=j$, and}\\
(1,&1,&1/2,1/2,\ldots)&\mbox{if $x>j$.}
\end{array}\right.\]
Now let $j\geq 1$. Then
\begin{eqnarray}
\lefteqn{P[T_{j+1}-T_j\geq j]}\nonumber\\
&\geq&\EE\left[P_{0,\om}\left[T_{j+1}-T_j\geq j,\ X_{T_{j-1}+1}=j-2\right], A_{j}\right]\nonumber\\
&=&\EE\left[E_{0,\om}\left[P_{j,\psi\left(\om,H_{T_j}\right)}
\left[T_{j+1}\geq j\right],\ X_{T_{j-1}+1}=j-2\right], A_{j}\right].\label{kata}
\end{eqnarray}\abel{kata}
Observe that 
\[P_{j,\psi\left(\om,H_{T_j}\right)}
\left[T_{j+1}\geq j\right]= P_{j,\om_j'}
\left[T_{j+1}\geq j\right],\]
where $\om_j'(x):=\psi\left(\om,H_{T_j}\right)(x)$ for $x\geq 0$ and 
$\om_j'(x)=\bar{\om}_j(x)$ for $x<0$. Moreover, $\om_j'\leq \bar{\om}_j$ 
on $\{\om(j,1)<1-\eps, X_{T_{j-1}+1}=j-2\}$.
Consequently, due to Lemma \ref{monoenv}, (\ref{kata}) is greater than or equal to
\[P_{j,\bar{\om}_j}[T_{j+1}\geq j] \ \EE\left[P_{0,\om}\left[X_{T_{j-1}+1}=j-2\right], A_j\right]
\ \geq \ P_{0,\bar{\om}_0}[T_{1}\geq j]\ \eps\ \al.\]
Hence it suffices to show that $\sum_{j\geq 1}P_{0,\bar{\om}_0}[T_1\geq j]=E_{0,\bar{\om}_0}[T_1]=\infty$.
Since
\[T_1=\sum_{k\geq 0}(T_{-k-1}\wedge T_1-T_{-k})\won\{T_{-k}<T_1\}\] 
we have
\begin{equation}\label{typo}
E_{0,\bar{\om}_0}[T_1]=\sum_{k\geq 0}E_{0,\bar{\om}_0}\left[
E_{0,\bar{\om}_0}\left[T_{-k-1}\wedge T_1-T_{-k}\mid\F_{T_{-k}}\right],\ T_{-k}<T_1\right].
\end{equation}\abel{typo}
For $k\geq 1$, the conditional expectation in (\ref{typo}) is on $\{T_{-k}<T_1\}$
by the strong Markov property equal to
\begin{eqnarray}
\lefteqn{
E_{-k,\psi\left(\bar{\om}_0,H_{T_{-k}}\right)}[T_{-k-1}\wedge T_1]
}\nonumber\\
&=&1+E_{-k+1,\psi\left(\bar{\om}_0,H_{T_{-k}+1}\right)}[T_{-k-1}\wedge T_1]\label{morning}\\
&\geq &E_{-k+1,\psi\left(\bar{\om}_0,H_{T_{-k}+1}\right)}[T_{-k-1}\wedge T_0]\ =\ 
2(k-1).\label{night}
\end{eqnarray}\abel{morning}\abel{night}
Here (\ref{morning}) holds because of $\bar{\om}_0(-k,1)=1$.  Equation (\ref{night}) is true
since the walker eats while traveling from 0 to $-k$ and back to $-k+1$ all the  cookies 
between $-k$ and $-1$, which have strength $>1/2$,  so that the  formula for the expected exit time of
a simple symmetric random walk from an interval (e.g.\ \cite[Ch.\ 14.3 (3.5)]{fell1})
can be applied.
Substituting this into (\ref{typo}) yields
\begin{eqnarray}
E_{0,\bar{\om}_0}[T_1]&\geq& 2\sum_{k\geq 1}(k-1)P_{0,\bar{\om}_0}[T_{-k}<T_1]\nonumber \\
&=& 2(1-\eps)
\sum_{k\geq 1}(k-1)P_{-1,\bar{\om}_0}[T_{-k}<T_0].\nonumber
\end{eqnarray}
Therefore, since the harmonic series diverges it suffices to show for the proof of $E_{0,\bar{\om}_0}[T_1]=\infty$ that for all $k\geq 2$,
\begin{equation}
P_{-1,\bar{\om}_0}[T_{-k}<T_0]=\frac{1}{(k-1)k}.\label{koko}
\end{equation}\abel{koko}
This is done by induction over $k$. For $k=2$, the left-hand side of (\ref{koko}) is $\bar{\om}_0(-1,1)=1/2$ by
definition of $\bar{\om}_0$. Now assume that (\ref{koko}) has been proven for $k$.
Then 
\begin{eqnarray*}
P_{-1,\bar{\om}_0}[T_{-k-1}<T_0]&=&E_{-1,\bar{\om}_0}\left[T_{-k}<T_0, 
P_{-k+1,\psi\left(\bar{\om}_0,H_{T_{-k}+1}\right)}
[T_{-k-1}<T_0]\right].
\end{eqnarray*}
As above, on $\{T_{-k}<T_0\}$ all the  cookies with strength $>1/2$
have been removed  by time $T_{-k}+1$ from the interval between $-k$ and $-1$.  
Therefore, the last expression equals
\[P_{-1,\bar{\om}_0}[T_{-k}<T_0]\frac{k-1}{k+1}
 = \frac{1}{(k-1)k}\ \frac{k-1}{k+1}=\frac{1}{k(k+1)}
\]
by induction hypothesis.
\end{proof}

The following example shows that the assumption
$\PP[\om(0,1)<1, \om(1,1)<1]>0$ of Theorem \ref{v0} is essential.

\begin{ex}\label{period}{\rm
Let $\om(x), x\in\Z,$ alternate between $(1,1,1/2,1/2,\ldots)$ and
 $(p,1,1/2,1/2,\ldots)$ where $1/2\leq p<1$ is fixed. Then $P_0$-a.s.\
$T_{k+1}-T_k\leq 3$ for all $k\geq 0$, which generates  a strictly
positive speed.}\hfill $\Box$
\end{ex}\abel{period}

\begin{problem} 
{\rm Of course, it is possible to generate a strictly positive speed $v$ by 
choosing $\om(x,i)\geq 1/2+\eps$ $\PP$-a.s.\ for all 
$x\in\Z^d, i\geq 1$, where $\eps>0$ is fixed. However, are a finite number of cookies with strength $>1/2$ (and $<1$)
per site
already sufficient to generate
positive speed? More precisely, for which integers $m\geq 3$, if any,
 is there some $1/2<p<1$ such that $v>0$ if $\PP$-a.s.\ for all $x\in\Z$, $\om(x,i)=p$ for $i\leq m$ and 
$\om(x,i)=1/2$ for $i>m$? 
}
\end{problem}



\bibliographystyle{amsalpha}
\vspace*{5mm}
{\sc \small
Department of Mathematics\\
Stanford University\\
Stanford, CA 94305, U.S.A.\\
E-Mail: {\rm zerner@math.stanford.edu} }
\end{document}